\documentclass{article} 

\usepackage{verbatim}
\usepackage{natbib}
\usepackage{amsfonts}
\usepackage{amssymb}
\usepackage{amsmath}
\usepackage{mathrsfs}
\usepackage{color}                 
\usepackage{graphicx}           
\usepackage{pst-func}  
\usepackage{pstricks-add} 
\usepackage{rotating}
\usepackage{setspace}
\usepackage{enumerate}
\usepackage[margin=3.5cm]{geometry}

\bibpunct{(}{)}{;}{a}{,}{,}

\newtheorem{as}{Assumption}
\newtheorem{lem}{Lemma}
\newtheorem{defi}{Definition}
\newtheorem{prop}{Proposition}

\newtheorem{theo}{Theorem}

\DeclareMathOperator{\sign}{sign}
 
\DeclareMathOperator*{\argmin}{argmin\;} 

\DeclareMathOperator{\Var}{Var}

\DeclareMathOperator{\diag}{diag}

\title{Uniformity and the delta method\thanks{This paper benefited greatly from the contributions of Jos\'e L. Montiel Olea as well as from discussions with Isaiah Andrews, Gary Chamberlain, Jinyong Hahn, Kei Hirano, Anna Mikusheva, James Stock, and Elie Tamer, and from the remarks of several anonymous referees.}}

\author{Maximilian Kasy\footnote{Assistant professor, Department of Economics, Harvard University.
Address: Littauer Center 200, 1805 Cambridge Street, Cambridge, MA 02138. email: maximiliankasy@fas.harvard.edu.}}

\begin{document}
\maketitle

\begin{abstract}
When are asymptotic approximations using the delta-method uniformly valid?
We provide sufficient conditions as well as closely related necessary conditions for uniform negligibility of the remainder of such approximations.
These conditions are easily verified and permit to identify settings and parameter regions where pointwise asymptotic approximations perform poorly.
Our framework allows for a unified and transparent discussion of uniformity issues in various sub-fields of econometrics. 
Our conditions involve uniform bounds on the remainder of a first-order approximation for the function of interest.
\end{abstract}

\noindent \textsc{JEL Classification:  C10\\
Key Words: Asymptotic theory, uniformity, delta method}

\clearpage

\section{Introduction}

Many econometric procedures are motivated and justified using asymptotic approximations.
Standard asymptotic theory provides approximations for fixed parameter values, letting the sample size go to infinity.
Procedures for estimation, testing, or the construction of confidence sets are considered justified if they perform well for large sample sizes, for any given parameter value.

Procedures that are justified in this sense might unfortunately still perform poorly for arbitrarily large samples.
This happens if the asymptotic approximations invoked are not \emph{uniformly} valid. In that case there are parameter values for every sample size such that the approximation is poor, even though for every given parameter value the approximation performs well for large enough sample sizes.
Which parameter values cause poor behavior might depend on sample size, so that poor behavior does not show up in standard asymptotics.
If a procedure is not uniformly valid, this can lead to various problems, including (i) large bias and mean squared error for estimators, (ii) undercoverage of confidence sets, and (iii) severe size distortions of tests.

Uniformity concerns are central to a  number of sub-fields of econometrics. 
The econometrics literature has mostly focused on uniform size control in testing and the construction of confidence sets.
Uniform validity of asymptotic approximations is however a more general issue, and is important even if we are not interested in uniform size control, but instead have decision theoretic or other criteria for the quality of an econometric procedure.
Literatures that have focused on uniformity issues include the literature on weak instruments, eg. \cite{SS97}, the literature on inference under partial identification in general and moment inequalities in particular, eg. \cite{imbens2004confidence}, and the literature on pre-testing and model selection, eg. \cite{leeb:2005}.

The purpose of this paper is to provide a unified perspective on failures of uniformity. 
We argue that in many settings the poor performance of estimators or lack of
uniformity of tests and confidence sets arises as a consequence of the lack of uniformity of approximations using the ``delta method.''
Motivated by this observation, we provide sufficient and necessary conditions for uniform negligibility of the remainder of an asymptotic approximation using the delta method. These conditions are easily verified.
This allows to spot potential problems with standard asymptotics and to identify parts of the parameter space where problems might be expected.  
Our sufficient conditions require that the function $\phi$ of interest is continuously differentiable, and that the remainder of the first order approximation $\phi(x + \Delta x) \approx \phi(x) + \phi'(x) \Delta x$ is uniformly small relative to the leading term $\phi'(x) \Delta x$, in a sense to be made precise below.

In the case of weak instruments, this condition fails in a neighborhood of $x_2 =0$ for the function $\phi(x) = x_1 / x_2$, which is applied to the ``reduced form'' covariances of instrument and outcome, and instrument and treatment.
In the case of moment inequalities or multiple hypothesis testing, remainders are not negligible in a neighborhood of kinks of the null-region,  where $\phi$ is for instance the distance of a statistic to the null-region.
For interval-identified objects as discussed by \cite{imbens2004confidence}, such a kink corresponds to the case of point-identification.
In the case of minimum-distance estimation, with over-identified parameters, remainders are non-negligible when the manifold traced by the model has kinks or high curvature.
In the case of pre-testing and model selection, this condition fails in a neighborhood of critical values for the pretest, where $\phi$ is the mapping from sample-moments to the estimated coefficients of interest; in the case of Lasso in the neighborhood of kinks of the mapping from sample-moments to the estimated coefficients.\footnote{Additional complications in pre-testing, model selection and Lasso settings arise because of drifting critical values or penalty parameters, so that they are only partially covered by our basic argument.}

The rest of this paper is structured as follows.
Section \ref{sec:litreview} provides a brief review of the literature.
Section \ref{sec:preliminaries} reviews definitions and discusses some preliminary results, including a result relating uniform convergence in distribution to uniformly valid confidence sets, and a uniform version of the continuous mapping theorem.
Section \ref{sec:uniformdelta} presents our central result, the sufficient and necessary conditions for uniform validity of the delta method. This section also shows that continuous differentiability on a compact domain is sufficient for uniform validity of the delta-method.
Section \ref{sec:applications} discusses several applications to illustrate the usefulness of our approach, including a number of stylized examples, weak instruments, moment inequalities, and minimum distance estimation.
Appendix \ref{sec:proofs} contains all proofs.

\section{Literature}

\label{sec:litreview}

Uniformity considerations have a long tradition in the statistical literature, at least since Hodges discussed his estimator and Wald analyzed minimax estimation.
Uniformity considerations have motivated the development of much of modern asymptotic theory and in particular the notion of limiting experiments, as reviewed for instance in \cite{le2012asymptotics}.

The interest in uniform asymptotics in econometrics was prompted by the poor finite-sample performance of some commonly used statistical procedures. Important examples include the study of `local-to-zero' behavior of estimators, tests and confidence sets in linear instrumental variable regression (\citealt{SS97, Moreira03, AMS06}); the `local-to-unit root' analysis for autoregressive parameters (see \citealt{stock1996confidence, mikusheva2007uniform});  and the behavior of estimators, tests and confidence sets that follow a pre-testing or a model selection stage (\citealt{leeb:2005, gugg10Haus}).

Much of this literature  has been concerned with finding statistical procedures that control size uniformly in large samples over a reasonable parameter space.
Our paper has a different objective. We argue that in most of these problems there are reduced-form statistics satisfying uniform central limit theorems (CLT) and uniform laws of large numbers (LLN).\footnote{Uniform CLT and uniform LLN in this paper are used to describe results that guarantee uniform convergence in distribution or probability of random vectors, rather than results that guarantee convergence of empirical processes.} The failure of some commonly used tests and confidence sets to be uniformly valid, despite the uniform convergence of reduced-form statistics, is a consequence of the lack of uniformity of the delta method.

There are some discussions of uniformity and the delta method in the literature.
\cite{van2000asymptotic}, for instance, discusses uniform validity of the delta-method in section 3.4. His result requires continuous differentiability of the function of interest $\phi$ on an open set, and convergence of the sequence of parameters $\theta_n$ to a point $\theta$ inside this open set. This result does not allow to study behavior near boundary points of the domain, which will be of key interest to us, as discussed below. The result in \cite{van2000asymptotic} section 3.4 is an implication of our more general theorems below.
Uniformity of the delta method has also been studied recently by \cite{belloni2013program} with a focus on infinite dimensional parameters. They provide a sufficient condition (a slight modification of the notion of uniform Hadamard differentiability in \cite{vandevaart1996weak}, p. 379) that guarantees the uniform validity of delta method approximations. 
In contrast to their condition, (i) we do not require the parameter space to be compact\footnote{Compactness excludes settings where problems arise near boundary points, such as weak instruments.}, and (ii) we provide necessary as well as sufficient conditions for uniform convergence.
The analysis in \cite{andrewsmikusheva2014} is closely related to ours in spirit. They consider tests of moment equality restrictions when the null-space is curved. First-order (delta-method) approximations to their test statistic are poor if the curvature is large.

The uniform delta-method established in this paper does not allow the function $\phi(\cdot)$ to depend on the sample size. \cite{Phillips2012} has extended the \emph{pointwise} delta-method in such a direction using an asymptotically locally relative equicontinuity condition.

Not all uniformity considerations fall into the framework discussed in our paper.
This is in particular true for local to unit root analysis in autoregressive models.
Model selection, and the Lasso, face problems closely related to those we discuss. Additional complications arise in these settings, however, because of drifting critical values or penalty parameters, which lead to ``oracle-property'' type approximations that are not uniformly valid.

\section{Preliminaries}

\label{sec:preliminaries}

In this section we introduce notation, define notions of convergence, and state some basic results (which appear to be known, but are scattered through the literature).
Throughout this paper, we consider random variables defined on the fixed sample space $\Omega$, which is equipped with a sigma-algebra $\mathscr{F}$, and a family of probability measures $P_\theta$ on $(\Omega, \mathscr{F})$ indexed by $\theta\in \Theta$.
$S,T, X, Y, Z$ and $S_n, T_n, X_n$ are random variables or random vectors defined on $\Omega$.
We are interested in asymptotic approximations with respect to $n$.
$\mu=\mu(\theta)$ denotes some finite-dimensional function of $\theta$, and $F$ is used to denote cumulative distribution functions, 
The derivative of $\phi(m)$ with respect to $m$ is denoted by $D(m)=\partial \phi / \partial m=\partial_m \phi$.

The goal of this paper is to provide conditions that guarantee uniform convergence in distribution.
There  are several equivalent ways to define convergence in distribution. One definition requires convergence in terms of the so called bounded Lipschitz metric, cf. \cite{vandevaart1996weak}, p73.
This definition is useful for our purposes, since it allows for a straightforward extension to uniform convergence.
\begin{defi}[Bounded Lipschitz metric]$\;$\\
Let $\mathbf{BL_1}$ be the set of all real-valued functions $h$ on $\mathbb{R}^{d_x}$ such that $ |h(x)| \leq 1$ and $|h(x) - h(x')| \leq \|x-x'\|$ for all $x,x'$.\\
The bounded Lipschitz metric on the set of random vectors with support in $\mathbb{R}^{d_x}$ is defined by
\begin{equation}
d^{\theta}_{BL} (X_1, X_2) := \sup_{h \in \mathbf{BL_1}} \left | E^{\theta}[h(X_1)] - E^{\theta}[h(X_2)] \right |.
\end{equation}
\end{defi}
In this definition, the distance of two random variables $X_1$ and $X_2$ depends on $\theta$, which indexes the distribution of $X_1$ and $X_2$, and thus also the expectation of functions of these random variables.
Standard asymptotics is about convergence for any \emph{given} $\theta$. Uniformity requires convergence for any \emph{sequence} of $\theta_n$, as in the following definition.

\begin{defi}[Uniform convergence]$\;$\\
\label{def:uniform}
Let $\theta \in \Theta$ be a (possibly infinite dimensional) parameter indexing the distribution of both $X_n$ and $Y_n$.
\begin{enumerate}
\item We say that $X_n$ converges uniformly in distribution to $Y_n$ if
\begin{equation}
d^{\theta_n}_{BL} (X_n, Y_n)  \rightarrow 0 
\end{equation}
for all sequences $\{\theta_n \in \Theta\}$.

\item We say that $X_n$ converges uniformly in probability to $Y_n$ if
\begin{equation}
 P^{\theta_n}(\|X_n - Y_n\| >\epsilon) \rightarrow 0 
\end{equation} 
for all $\epsilon > 0$ and all sequences $\{\theta_n \in \Theta\}$.
\end{enumerate}
\end{defi}

\textbf{Remarks:}
\begin{itemize}
\item As shown in \cite[][section 1.12]{vandevaart1996weak}, convergence in distribution of $X_n$ to $X$, defined in the conventional way as convergence of cumulative distribution functions (CDFs) at points of continuity of the limiting CDF, is equivalent to convergence of $d^{\theta}_{BL} (X_n, X)$ to $0$.
\item Our definition of uniform convergence might seem slightly unusual. We define it as convergence of a sequence $X_n$ toward another sequence $Y_n$. In the special case where $Y_n=X$ so that $Y_n$ does not depend on $n$, this definition reduces to the more conventional one, so our definition is more general.
\item There are several equivalent ways to define uniform convergence, whether in distribution or in probability. 
Definition \ref{def:uniform} requires convergence along all sequences $\theta_n$.
The following Lemma \ref{lem:characterization}, which is easy to prove, shows that this is equivalent to requiring convergence of suprema over all $\theta$.
\end{itemize}

\begin{lem}[Characterization of uniform convergence]$\;$
\label{lem:characterization}
\begin{enumerate}
\item $X_n$ converges uniformly in distribution to $Y_n$ if and only if
\begin{equation}
 \sup_{\theta\in \Theta} d^\theta_{BL} (X_n, Y_n) \rightarrow 0 
\end{equation}

\item $X_n$ converges uniformly in probability to $Y_n$ if and only if
\begin{equation}
\sup_{\theta\in \Theta} P^\theta(\|X_n -Y_n\| >\epsilon) \rightarrow 0
\end{equation}
for all $\epsilon > 0$.

\end{enumerate}

\end{lem} 
 The proof of this lemma and of all following results can be found in appendix \ref{sec:proofs}.\\


%
%

Uniform convergence safeguards, in large samples, against asymptotic approximations performing poorly for some values of $\theta$. Absent uniform convergence, there are $\theta_n$ for arbitrarily large $n$ for which the approximation is far from the truth.
Guaranteeing uniformity is relevant, in particular, to guarantee the validity of inference procedures.
The following result shows how uniform convergence of a test-statistic to a pivotal distribution allows to construct confidence sets with appropriate coverage.
This result could equivalently be stated in terms of hypothesis testing.

 \begin{lem}[Uniformly valid confidence sets]$\;$\\
 \label{lem:uniforminference}
Suppose that $Z_n$ converges uniformly in distribution to $Z$, where $Z_n = Z_n(\mu)$.
Suppose further that $Z$ is continuously distributed and pivotal, that is, the distribution of $Z$ does not depend on $\theta$.
Let $z$ be the $1-\alpha$ quantile of the distribution of $Z$.
Then
\begin{equation}
C_n := \{m:  Z_n(m) \leq z\}
\end{equation}
is such that
\begin{equation}
P^{\theta_n}(\mu(\theta_n) \in C_n) \rightarrow 1-\alpha
\end{equation}
for any sequence $\theta_n\in \Theta$.
\end{lem}

Lemma \ref{lem:uniforminference} establishes the connection between our definition of uniform convergence in distribution and uniformly valid inference. The latter hinges on convergence of 
\[\sup_{\theta\in \Theta}  \left | F^\theta_{Z_n(m)}(z)-F_{Z}(z) \right |\] to $0$ for a given critical value $z$, whereas uniform convergence in distribution of $Z_n$ to $Z$ can be shown to be equivalent to convergence of
\[\sup_{\theta\in \Theta} \sup_z \left | F^\theta_{Z_n(m)}(z)-F_{Z}(z) \right |\]
to $0$. If we were to require uniform validity of inference for arbitrary critical values, this would in fact be equivalent to uniform convergence in distribution of the test-statistic.
We should emphasize again, however, that uniform convergence in distribution is a concern even if we are not interested in uniform size control, but for instance in the risk of an estimator.\\

From our definition of uniform convergence, it is straightforward to show the following uniform version of the continuous mapping theorem.
The standard continuous mapping theorem states that convergence in distribution (probability) of $X_n$ to $X$ implies convergence in distribution (probability) of $\psi(X_n)$ to $\psi(X)$ for any \emph{continuous} function $\psi$.
Our uniform version of this result needs to impose the slightly stronger requirement that $\psi$ be uniformly continuous (for uniform convergence in probability) or \emph{Lipschitz-continuous} (for uniform convergence in distribution).
\begin{theo}[Uniform continuous mapping theorem]$\;$\\
\label{theo:uniformCMT}
Let $\psi(x)$ be a function of $x$ taking values in $\mathbb{R}^l$.
\begin{enumerate}
\item Suppose $X_n $ converges uniformly in distribution to $Y_n$.\\
If $\psi$ is Lipschitz-continuous, then $\psi(X_n)$ converges uniformly in distribution to $\psi(Y_n)$.
\item Suppose $X_n $ converges uniformly in probability to $Y_n$.\\
If $\psi$ is uniformly continuous, then $\psi(X_n)$ converges uniformly in probability to $\psi(Y_n)$. 
\end{enumerate}
\end{theo}

\textbf{Remarks:}
\begin{itemize}
\item 
Continuity of $\psi$ would not be enough for either statement to hold.
To see this, consider the following example: assume $\theta \in \mathbb{R}^+$, $X_n = \theta$, and $Y_n = X_n + 1/n$. Then clearly $Y_n$ converges uniformly (as a sequence, in probability, and in distribution) to $X_n$.
Let now $\psi(x) = 1/x$, and $\theta_n = 1/n$. $\psi$ is a continuous function on the support of $X_n$ and $Y_n$. Then $\psi(X_n)=1$, $\psi(Y_n)=1/2$, and $P^{\theta_n}(|\psi(X_n) - \psi(Y_n)| = 1/2) = 1$, and thus $\psi(Y_n)$ does not converge uniformly (in probability, or in distribution) to $\psi(X_n)$.

\item There is, however, an important special case for which continuity of $\psi$ \emph{is} enough: If $Y_n = Y$ and the distribution of $Y$ does not depend on $\theta$, so that $Y$ is a pivot, then convergence of $X_n$ to $Y$ uniformly in distribution implies that $\psi(X_n)$ converges uniformly in distribution to $\psi(Y)$ for any \emph{continuous} function $\psi$.
This follows immediately from the standard continuous mapping theorem, applied along arbitrary sequences $\{\theta_n\}$.
\end{itemize}

\section{The uniform delta-method}

\label{sec:uniformdelta}

We will now discuss the main result of this paper.
In the following theorem \ref{theo:uniformdelta}, we consider a sequence of random variables (or random vectors) $T_n$  such that 
\[S_n:=r_n (T_n - \mu) \rightarrow^d S\]
 uniformly.
We are interested in the distribution of some function $\phi$ of $T_n$. 
Let
\begin{align*}
D(m) &=  \frac{\partial \phi}{\partial m}(m), \textrm{ and}\\
E(m) &= \diag(\|D_k(m)\|)^{-1},
\end{align*} 
where $\|D_k(m)\|$ is the norm of the $k$th row of $D(m)$.
Consider the normalized sequence of random variables
\begin{equation}
X_n = r_n E(\mu) (\phi(T_n) - \phi(\mu)).
\end{equation}
We aim to approximate the distribution of $X_n$ by the distribution of  
\begin{equation}
X:=  E(\mu) D(\mu)\cdot S
\end{equation}
Recall that the distributions of $T_n$ and $X_n$ are functions of both $\theta$ and $n$.
$\mu$ and the distribution of $S$ are functions of $\theta$ (cf. definition \ref{def:uniform}). The sequence $r_n$ is not allowed to depend on $\theta$; the leading case is $r_n = \sqrt{n}$.
If $T_n$ is a vector of dimension $d_t$ and $X_n$ is of dimension $d_x$, then $D=\frac{\partial \phi}{\partial \mu}$ is a $d_x \times d_t$ matrix of partial derivatives.
$E(\mu) = \diag(\|D_k(\mu)\|)^{-1}$ is a $d_x \times d_x$ diagonal matrix which serves to normalize the rows of $D(\mu)$.\\

Our main result, theorem \ref{theo:uniformdelta}, requires that the ``reduced form'' statistics $S_n$ converge uniformly in distribution to a tight family of continuously distributed random variables $S$.
Uniform convergence of the reduced form can be established for instance using central limit theorems for triangular arrays, cf. lemma \ref{lem:uniformCLT} below.
\begin{as}[Uniform convergence of $S_n$]$\;$\\
\label{as:uniformS}
Let $S_n:=r_n (T_n - \mu)$.
\begin{enumerate}
\item $S_n \rightarrow S$ uniformly in distribution for a sequence $r_n \rightarrow \infty$ which does not depend on $\theta$.

\item $S$ is continuously distributed for all $\theta \in \Theta$.

\item  The collection $\{S(\theta)\}_{\theta \in \Theta}$ is tight, that is, for all $\epsilon>0$ there exists an $M < \infty$ such that $P(\|S\| \leq M) \geq 1-\epsilon$ for all $\theta$.
 
\end{enumerate}
\end{as}

The leading example of such a limiting distribution is the normal distribution $S\sim N(0, \Sigma(\theta))$. This satisfies the condition of tightness if there is a uniform upper bound on $ \|\Sigma(\theta)\|$.\\

The sufficient condition for uniform convergence in theorem \ref{theo:uniformdelta} below furthermore requires a uniform bound on the remainder of a first order approximation to $\phi$.
Denote the normalized remainder of a first order Taylor expansion of $ \phi $ around $m$ by
\begin{equation} 
\Delta(t, m) = \frac{1}{\|t - m\|}\left \|E(m)\cdot(\phi(t)- \phi(m) -  D(m) \cdot(t - m)) \right \|.
\end{equation}
The function $\phi$ is differentiable -- as required by the pointwise delta-method -- if and only if $\Delta(t, m)$ goes to $0$ as $t\rightarrow m$ for fixed $m$.

Note the role of $E(m)$ in the definition of $\Delta$. Normalization by $E(m)$ ensures that we are considering the magnitude of the remainder \emph{relative} to the leading term $D(m) \cdot(t - m)$. This allows us to consider settings with unbounded derivatives $D(m)$.\\

The first part of theorem \ref{theo:uniformdelta}  states a sufficient condition for uniform convergence of $X_n$ to $X$. This condition is a form of ``uniform differentiability;'' it requires that the remainder $\Delta(t, m)$ of a first order approximation to $\phi$ becomes uniformly small relative to the leading term as $\|t-m\|$ becomes small.
This condition fails to hold in all motivating examples mentioned at the beginning of this paper: weak instruments, moment inequalities, model selection, and the Lasso.

The second part of theorem \ref{theo:uniformdelta}  states a condition which implies that   $X_n$ does not converge uniformly to $X$  in distribution. This condition requires the existence of a $m_n \in \mu(\Theta)$ such that the remainder of a first order approximation becomes large relative to the leading term.


\begin{theo}[Uniform delta method]$\;$\\
\label{theo:uniformdelta}
Suppose assumption \ref{as:uniformS} holds.
 Let $\phi$  be a function which is continuously differentiable everywhere in an open set containing $\mu(\Theta)$, the set of $\mu$ corresponding to the parameter space $\Theta$.\footnote{Note that $\Theta$ and  $\mu(\Theta)$ might be open and/or unbounded.}
Assume that $D(\mu)$ has full row rank for all $\mu \in \mu(\Theta)$. 

\begin{enumerate}
\item Suppose that
\begin{equation}
\label{eq:vanishingremainder}
\Delta(t, m) \leq \widetilde{\delta}(\| t - m\|).
\end{equation} 
for some function $\widetilde{\delta}$ where $\lim_{\epsilon\rightarrow 0} \widetilde{\delta}(\epsilon) =0$, and for all $m\in \mu(\Theta)$.

Then $X_n$ converges uniformly in distribution to $X$.

\item Suppose there exists an open set $A\subset \mathbb{R}^{d_t}$ such that $\inf_{s \in \overline{A}} \|s\| = \underline{s} >0$ and $P_\theta(S\in A)  \geq \underline{p} >0$ for all $\theta\in \Theta$, and a sequence  $(\epsilon'_n, m_n)$, $\epsilon'_n>0$ and $m_n\in \mu(\Theta)$, such that
\begin{align}
\Delta(m_n + s / r_n, m_n)  &\geq \epsilon'_n \quad \forall s  \in A, \label{eq:divergingremainder}\\
\epsilon'_n & \rightarrow \infty. \nonumber
\end{align}
 
Then $X_n$ does not converge uniformly in distribution to $X$.

\end{enumerate}

\end{theo}


In any given application, we can check uniform validity of the delta method by verifying whether the sufficient condition in part 1 of this theorem holds. 
If it does not, it can in general be expected that uniform convergence in distribution will fail.
Part 2 of the theorem allows to show this directly, by finding a sequence of parameter values such that the remainder of a first-order approximation dominates the leading term.
There is also an intermediate case, where the remainder is of the same order as the leading term along some sequence $\theta_n$, so that the condition in part 2 holds for some constant rather than diverging sequence $\epsilon'_n$. This intermediate case is covered by neither part of theorem \ref{theo:uniformdelta}. In this intermediate case, uniform convergence in distribution would be an unlikely coincidence. Non-convergence for such intermediate cases is best shown on a case-by-case basis.
Section \ref{ssec:stylizedexamples} discusses several simple examples of functions $\phi(t)$ for which we demonstrate that the uniform validity of the delta method fails: $1/t$, $t^2$,  $|t|$, and $\sqrt{t}$.

The following theorem \ref{theo:sufficient}, which is a consequence of theorem \ref{theo:uniformdelta}, shows that a compact domain of $\phi$ is sufficient for condition \eqref{eq:vanishingremainder} to hold.
While compactness 
is too restrictive for most settings of interest, this result indicates where we might expect uniformity to be violated: either in the neighborhood of boundary points of the domain $\mu(\Theta)$ of $\phi$, if this domain is not closed, or as $m \rightarrow \infty$. The applications discussed below are all in the former category, and so are the functions  $1/t$, $t^2$, $|t|$,and $\sqrt{t}$.
Define the set $\mathscr{T}^\epsilon $, for any given set $\mathscr{T}$, as
\begin{align}
\mathscr{T}^\epsilon = \{ t&: \|t-\mu\| \leq \epsilon, 
\; \mu \in \mathscr{T}\}.
\end{align} 
\begin{theo}[Sufficient condition]$\;$\\
\label{theo:sufficient}
Suppose that $\mu(\Theta)$ is compact and  $\phi$ is everywhere continuously differentiable on $\mu(\Theta)^\epsilon$ for some $\epsilon>0$, and that $D(\mu)$ has full row rank for all $\mu \in \mu(\Theta)$.
Then condition \eqref{eq:vanishingremainder} holds.
\end{theo}

Theorem \ref{theo:uniformdelta} requires that the ``reduced form'' statistics $S_n$ converge uniformly in distribution to a tight family of continuously distributed random variables $S$.
One way to establish uniform convergence of reduced form statistics is via central limit theorems for triangular arrays, as in the following lemma, which immediately follows from Lyapunov's central limit theorem.
\begin{lem}[Uniform central limit theorem]$\;$\\
\label{lem:uniformCLT}
Let $Y_i$ be i.i.d. random variables with mean $\mu(\theta)$ and variance $\Sigma(\theta)$.
Assume that $E\left [\|Y_i^{2+\epsilon}\|\right ] <M$ for a constant $M$ independent of $\theta$. Then 
\[S_n :=\frac{1}{\sqrt{n}} \sum_{i=1}^n  \left (Y_i - \mu(\theta)\right )\]
converges uniformly in distribution to the tight family of continuously distributed random variables $S \sim N(0, \Sigma(\theta))$.
\end{lem}

In lemma \ref{lem:uniforminference} above we established that uniform convergence to a continuous pivotal distribution allows to construct uniformly valid hypothesis tests and confidence sets. Theorem \ref{theo:uniformdelta} guarantees uniform convergence in distribution; some additional conditions are required to allow for construction of a statistic which uniformly converges to a \emph{pivotal} distribution. The following proposition provides an example.  
\begin{prop}[Convergence to a pivotal statistic]$\;$\\
\label{prop:normalpivot}
Suppose that the assumptions of theorem \ref{theo:uniformdelta} and condition \eqref{eq:vanishingremainder} hold.
Assume additionally that 
\begin{enumerate}
\item  $\frac{\partial \phi}{\partial m}$ is Lipschitz continuous and its determinant is bounded away from $0$,
\item $S \sim N(0, \Sigma)$, where $\Sigma$ might depend on $\theta$ and its determinant is bounded away from $0$,
\item and that  $\widehat{\Sigma}$ is a uniformly consistent estimator for $\Sigma$,  in the sense that $\| \widehat{\Sigma} - \Sigma \|$ converges to $0$ in probability along any sequence $\theta_n$.
\end{enumerate}
Let
\begin{equation}
Z_n = \left (\frac{\partial \phi}{\partial t} (T_n) \cdot   \widehat{\Sigma} \cdot  \frac{\partial \phi}{\partial t} (T_n)'  \right  )^{-1/2} \cdot X_n.
\end{equation}
Then
\[Z_n \rightarrow N(0, I)\]
uniformly in distribution.

\end{prop}

%
%

\section{Applications}

\label{sec:applications}

\subsection{Stylized examples}
\label{ssec:stylizedexamples}

The following examples illustrate various ways in which functions $\phi(t)$ might violate the necessary and sufficient condition in theorem \ref{theo:uniformdelta}, and the sufficient condition of theorem \ref{theo:sufficient} (continuous differentiability on an $\epsilon$-enlargement of a compact domain).
In all of the examples we consider, problems arise near the boundary point $0$ of the domain $\mu(\Theta)=\mathbb{R} \setminus \{0\}$ of $\phi$.\footnote{Recall that we consider functions $\phi$ which are continuously differentiable on the domain $\mu(\Theta)$; $\phi$ might not be differentiable or not even well defined at boundary points of this domain.}
All of these functions might reasonably arise in various statistical contexts.

\begin{itemize}
\item
The first example, $1/t$, is a stylized version of the problems arising in inference using weak instruments.
This function diverges at $0$, and problems emerge in a neighborhood of this point.

\item
The second example, $t^2$, is seemingly very well behaved, and in particular continuously differentiable on all of $\mathbb{R}$. 
In a neighborhood of $0$, however, the leading term of a first-order expansion becomes small relative to the remainder term.

 \item The third example, $|t|$, is a stylized version of the problems arising in inference based on moment inequalities. This function is continuous everywhere on $\mathbb{R}$. It is not, however, differentiable at $0$, and problems emerge in a neighborhood of this point.

\item
The last example, $\sqrt{t}$ on $\mathbb{R}^+$, illustrates an intermediate case between  example 1 and example 3. This function is continuous on its domain and differentiable everywhere except at $0$; in contrast to $|t|$ it does not posses a directional derivative at $0$.
 \end{itemize}

In each of these examples problems emerge at a boundary point of the domain of continuous differentiability; such problems could not arise if the domain of $\phi$ were closed. The neighborhood of such boundary points is often of practical relevance.
Problems could in principle also emerge for very large $m$; such problems could not arise if the domain of $\phi$ were bounded. Very large $m$ might however be considered to have lower ``prior likelihood'' in many settings.
For each of these examples we provide analytical expressions for $\Delta$, a discussion in the context of theorem \ref{theo:uniformdelta}, as well as visual representations of $\Delta$.

\subsubsection{$\mathbf{\phi(t)= 1/t}$}

\begin{figure}[p] 
\centering
\caption{$\Delta(t,m)$ for $\phi(t)=1/t$.}
\includegraphics[keepaspectratio,  scale=1]{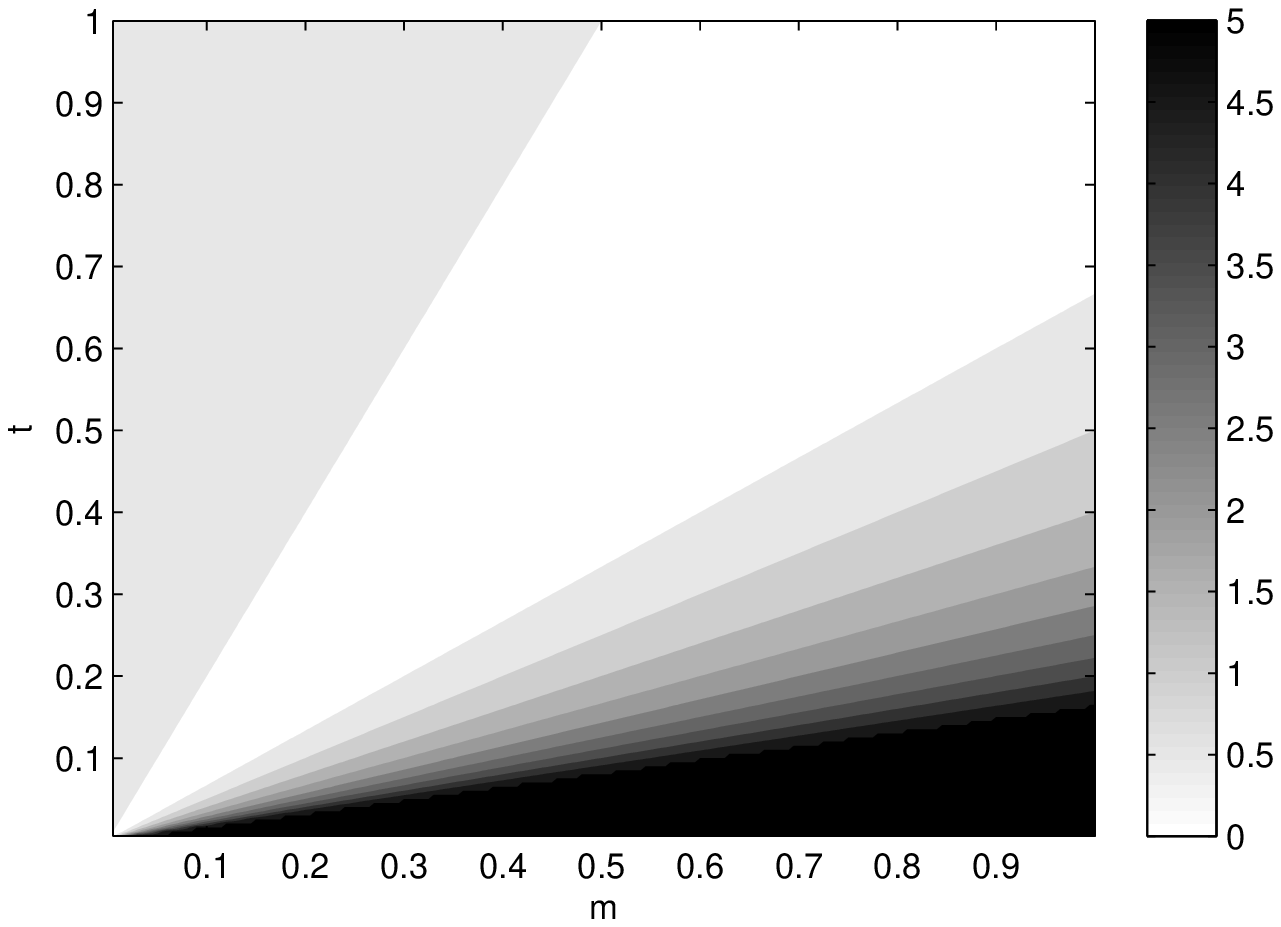}\\
\label{figure:Delta-1overt}
 
\centering
\caption{$\Delta(t,m)$ for $\phi(t)=t^2$.}
\includegraphics[keepaspectratio,  scale=1]{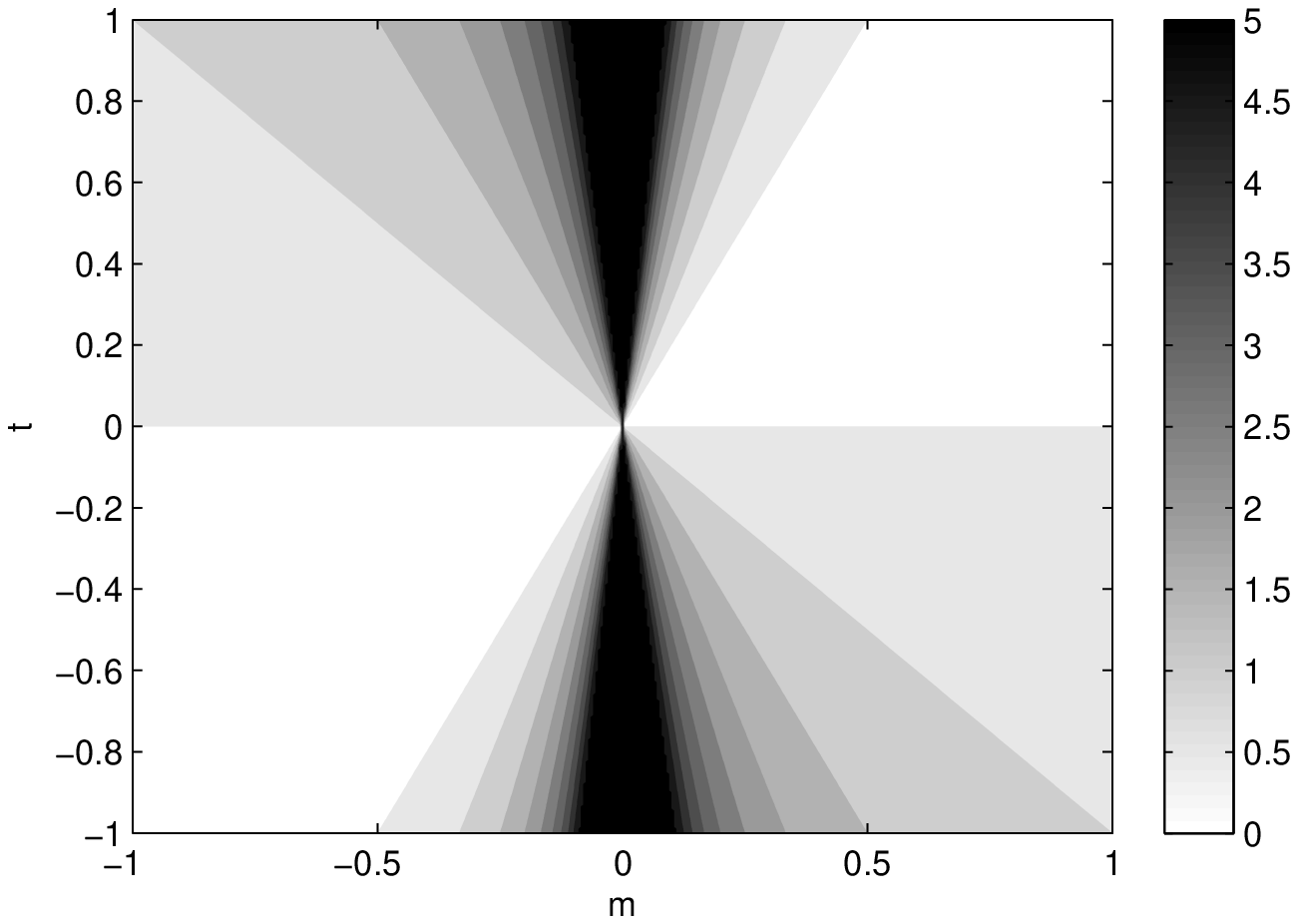}\\
\label{figure:Delta-square}
 
\end{figure}

 For the function $\phi(t)= 1/t$ we get $D(m)=\partial_m \phi(m) = - 1 / m^2$, $E(m)=m^2$, and
\begin{align*}
\Delta(t, m) &= \frac{m^2}{|t - m|}\left |  \frac{1}{t} - \frac{1}{m} + \frac{t-m}{m^2} \right | \\
&=\frac{m^2}{|t - m|}\left |\frac{ m\cdot (m- t) +   t\cdot(t-m)}{m^2 \cdot t} \right | \\
&= \left|  \frac{t-m}{ t} \right|.
\end{align*}  
Figure \ref{figure:Delta-1overt} shows a plot of the remainder $\Delta$; we will have similar plots for the following examples.

We get $\Delta > \epsilon'$ if and only if $|t-m| / |t| >\epsilon'$.
This holds if either
\begin{align*}
t &< \frac{m}{1+\epsilon'} \textrm{, or}\\
t &> \frac{m}{1- \epsilon'} \textrm{ and }   \epsilon' <1.
\end{align*}

We can show failure of the delta method to be uniformly valid for $\phi(t) =1/t$, using the second part of theorem \ref{theo:uniformdelta}.
In the notation of this theorem, let
\begin{align*}
r_n &= \sqrt{n}\\
\epsilon'_n &= \sqrt{n}\\
m_n &= \frac{1}{\sqrt{n}} + \frac{1}{n}\\
A&= ]-2, -1[.
\end{align*}
It is easy to check that the condition of theorem \ref{theo:uniformdelta}, part 2 applies for these choices.

\subsubsection{$\mathbf{\phi(t)= t^2}$}

For the function $\phi(t)= t^2$ we get $D(m)=\partial_m \phi(m) = 2 m$, $E(m)=1/(2 m)$, and
\begin{align*}
\Delta(t, m) &= \frac{1}{2m \cdot |t - m|}\left | t^2- m^2 - 2m \cdot (t-m)  \right | \\
&= \frac{1}{2m \cdot |t - m|}\left |(t-m)^2  \right | \\
&=\left |\frac{t-m}{2m}\right |
\end{align*}  
We therefore get $\Delta > \epsilon'$ if and only if $|t-m| / |2m| >\epsilon'$.
This holds if either
\begin{align*}
t &< m \cdot (1 - 2\epsilon')\textrm{, or}\\
t &> m \cdot (1 + 2\epsilon').
\end{align*}

We can again show failure of the delta method to be uniformly valid for $\phi(t) =t^2$, using the second part of theorem \ref{theo:uniformdelta}.
In the notation of this theorem, let
\begin{align*}
r_n &= \sqrt{n}\\
\epsilon'_n &= \sqrt{n} /2\\
m_n &= \frac{1}{n}\\
A&= ]1,2[.
\end{align*}
It is easy to check that the condition of theorem \ref{theo:uniformdelta}, part 2 applies for these choices.

\subsubsection{$\mathbf{\phi(t)= |t|}$}

 For the function $\phi(t)= |t|$, we get $D(m)=\partial_m\phi(m) = \sign(m)$, $E(m)=1$ and thus the normalized
remainder of the first-order Taylor expansion is given by
\begin{align*}
\Delta(t, m) &= \frac{1}{|t - m|}\left ||t| -|m| -  \sign(m)  \cdot(t - m) \right |\\
& =\mathbf{1}(t\cdot m \leq 0) \frac{2\cdot  |t|}{|t-m|}.
\end{align*} 

To see that the sufficient condition of the first part of theorem \ref{theo:sufficient} does not hold for this example, consider the sequence
\begin{align*}
m_n&= 1/n\\
t_n&= -1/n\\
\Delta(t_n, m_n) &= 1.
\end{align*}
 
In this example, however, the remainder does not diverge; the condition for non-convergence of the second part of  theorem \ref{theo:sufficient} does not apply either. To see this, note that
\[\Delta(t,m) \leq 2\]
 for all $t,m$ in this case.
 We are thus in the intermediate case, where remainder and leading term remain of the same order as $m$ approaches $0$, and have to show non-convergence ``by hand.''
 To do so, suppose $S_n$ converges uniformly in distribution to $S \sim N(0,1)$, and $r_n=\sqrt{n}$.
It immediately follows that $X\sim N(0,1)$ for all $\theta$. 
 Consider a sequence $\theta_n$ such that $m_n = \mu(\theta_n) =1/n$.
 For this sequence we get
 \begin{align*}
 F^{\theta_n}_{X_n} (0) &\rightarrow 0\\
F^{\theta_n}_{X} (0) &= 1/2,
 \end{align*}
 which immediately implies that we cannot have $X_n \rightarrow^d X$ along this sequence.

\subsubsection{$\mathbf{\phi(t)= \sqrt{t}}$}

For the function $\phi(t)= \sqrt{t}$, considered to be a function on $\mathbb{R}^+$,
 we get $D(m)=\partial_m \phi(m) = m^{-1/2} / 2$, $E(m)=2\cdot m^{1/2}$, and
\begin{align*}
\Delta(t, m) &=  \frac{2\cdot \sqrt{m}}{|(\sqrt{t} - \sqrt{m})(\sqrt{t} + \sqrt{m})|}\left |  \sqrt{t} - \sqrt{m} - \frac{t-m}{2 \sqrt{m}} \right | \\
&= \left | \frac{2\cdot \sqrt{m}}{\sqrt{t} + \sqrt{m}}-1 \right | \\
&= \left |  \frac{\sqrt{m}- \sqrt{t}}{\sqrt{m} +\sqrt{t} } \right |.
\end{align*}

This implies   $\Delta > \epsilon'$ if and only if $ |\sqrt{m}- \sqrt{t}| >\epsilon'  (\sqrt{m} +\sqrt{t}) $.
This holds if either
\begin{align*}
\sqrt{t} &<  \sqrt{m} \cdot  \frac{1- \epsilon' }{1+\epsilon'}  \textrm{, or}\\
\sqrt{t} &> \sqrt{m} \cdot  \frac{1+ \epsilon'}{1-\epsilon'}  \textrm{ and }  \epsilon' <1.
\end{align*}

Again, the sufficient condition of the first part of theorem \ref{theo:sufficient} does not hold for this example.
To show this, consider the sequence
\begin{align*}
m_n&= 4/n^2\\
t_n&= 1/n^2\\
\Delta(t_n, m_n) &= 1/3.
\end{align*}
Again, as well, the remainder does not diverge; the condition for non-convergence of the second part of  theorem \ref{theo:sufficient} does not apply either. To see this, note that
\[\Delta(t,m) \leq 1\]
 for all $t,m>0$.

To show non-convergence ``by hand,'' suppose $S_n$ converges uniformly in distribution to $S \sim \chi^2_1$, and $r_n=\sqrt{n}$.
It immediately follows that $X\sim \chi^2_1$ for all $\theta$. 
 Consider a sequence $\theta_n$ such that $m_n = \mu(\theta_n) =1/n$.
 For this sequence we get
 \begin{align*}
X_n & \rightarrow^{d,\theta_n}  N(0,1)\\
X &\sim^{\theta_n} \chi^2_1,
 \end{align*}
 which immediately implies that we cannot have $X_n \rightarrow^d X$ along this sequence.

\begin{figure}[p] 
\centering
\caption{$\Delta(t,m)$ for $\phi(t)=|t|$.}
\includegraphics[keepaspectratio,  scale=1]{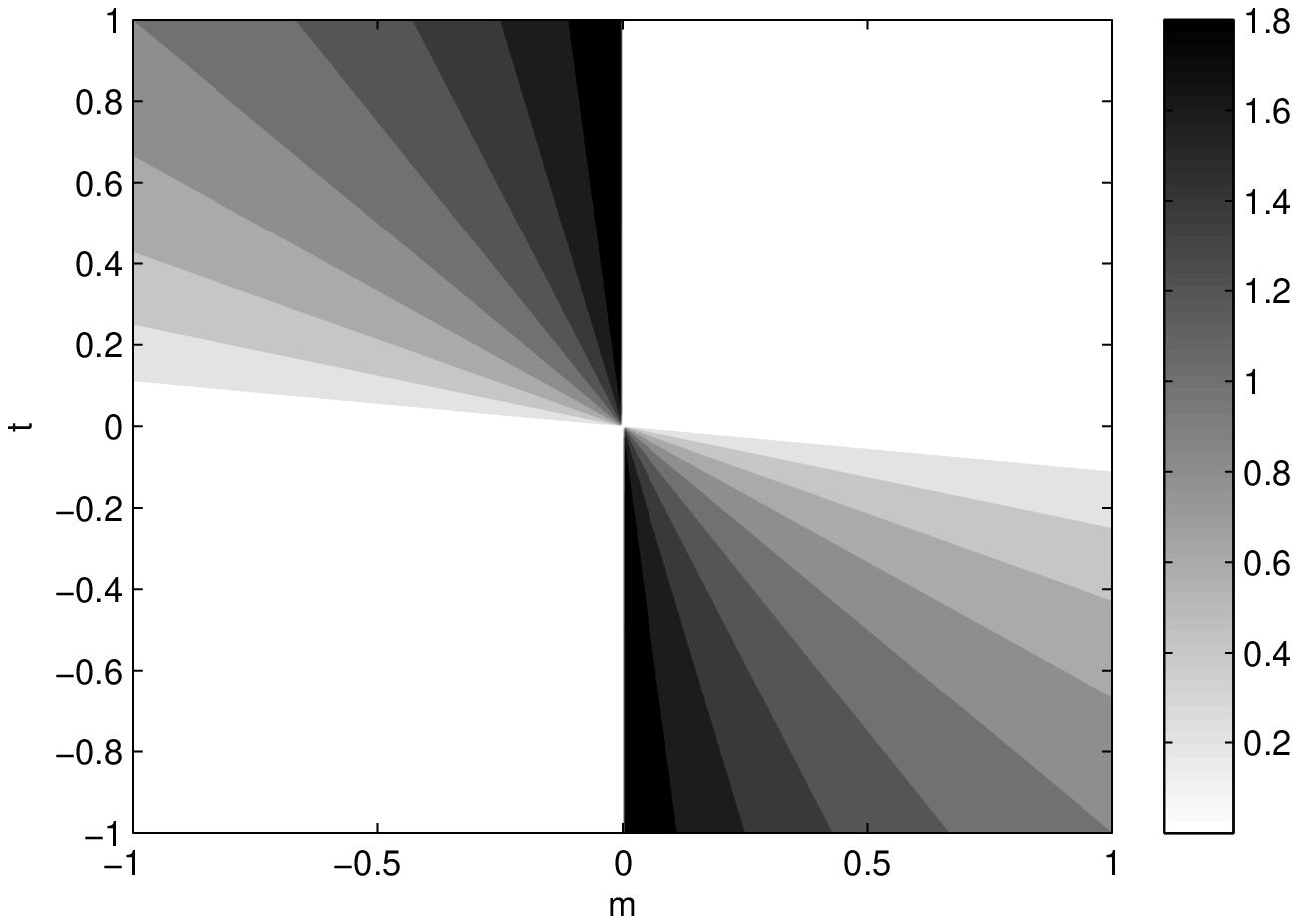}\\
\label{figure:Delta-absval}

\centering
\caption{$\Delta(t,m)$ for $\phi(t)=\sqrt{t}$.}
\includegraphics[keepaspectratio,  scale=1]{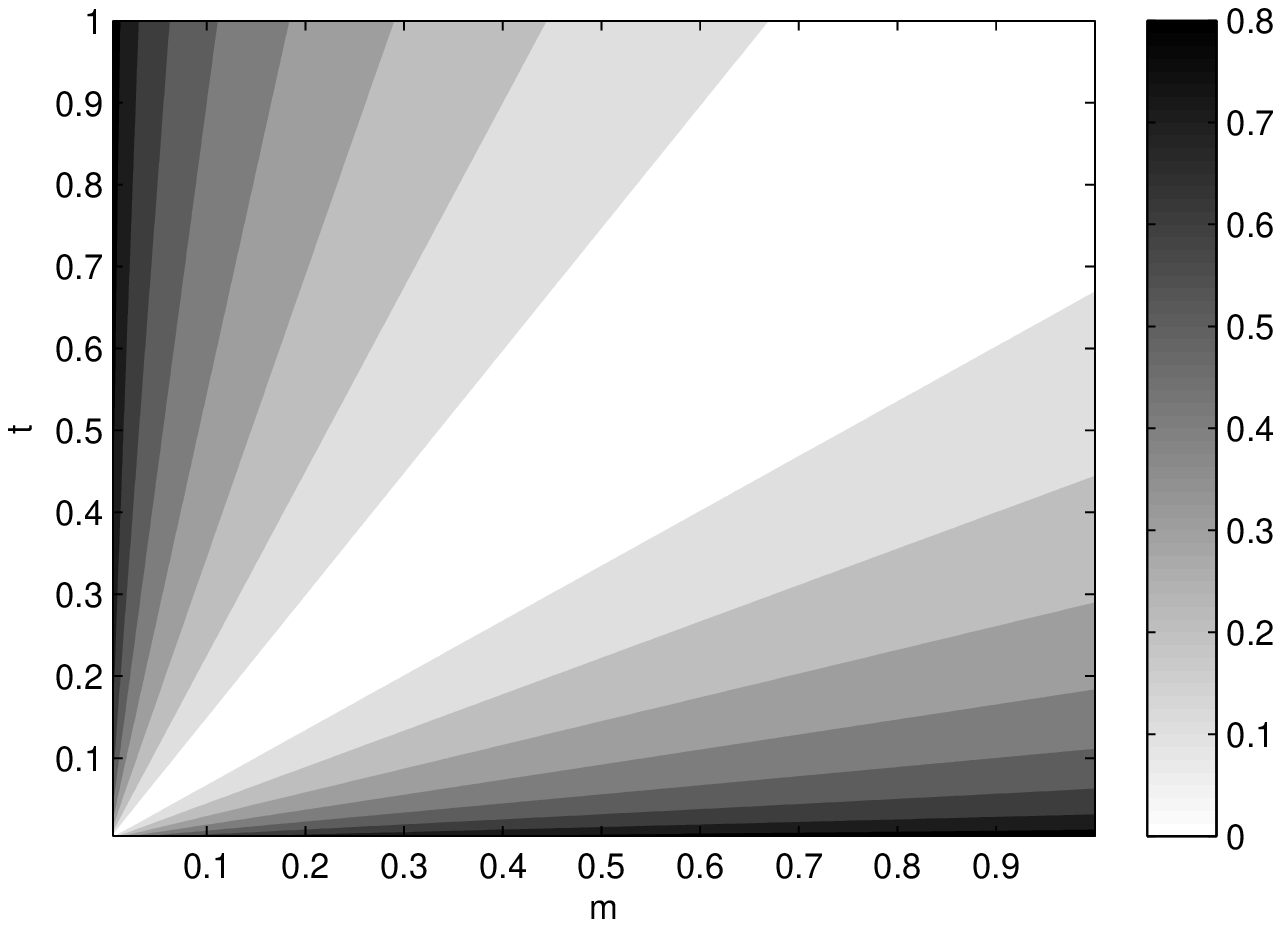}\\
\label{figure:Delta-sqrt}

\end{figure}

\subsection{Weak instruments}

Suppose we observe an i.i.d. sample $(Y_i, D_i, Z_i)$,
where we assume for simplicity that $E[Z] = 0$ and $\Var(Z)=1$.
Consider the  linear IV estimator
\begin{equation}
\widehat{\beta} :=  \frac{E_n[Z Y] }{ E_n[Z D]}.
\end{equation}
To map this setting into the general framework discussed in this paper, let
\begin{align}
T_n &= E_n[(Z Y, Z D)] \nonumber\\
\mu(\theta) &=  E[(Z Y, Z D)] \nonumber\\
\Sigma(\theta) &= \Var((ZY, DY)) \textrm{ and} \nonumber\\
\phi(t)&= \frac{t_1}{t_2}.
\end{align}
 In this notation, $\widehat{\beta} = \phi(T_n)$.
  This is a version of the ``weak instrument'' setting originally considered by
   \cite{SS97}.
Application of lemma \ref{lem:uniformCLT} to the statistic $T_n$ yields
\begin{equation*}
\sqrt{n} \cdot (T_n - \mu(\theta)) \rightarrow N(0, \Sigma(\theta)),
\end{equation*}
as long as $E[\|(ZY, DY)\|^{2+\epsilon}]$ is uniformly bounded  for some $\epsilon >0$.\\

Theorem \ref{theo:uniformdelta} thus applies. 
Taking derivatives we get $D(m) = (1/m_2, - m^1/m_2^2)$, and the inverse norm of $D(m)$ is given by $E(m) = \|D(m)\|^{-1} = m_2^2 / \|m\|$.
Some algebra, which can be found in appendix \ref{sec:proofs}, yields 
\begin{equation}
\Delta(t, m) =  \frac{1}{\|m\| \cdot \|t-m\|} \left | \frac{t_2-m_2}{t_2} \right | \: \:\left | m_2 \cdot (t_1-m_1) - m_1 \cdot (t_2-m_2)\right |.
\label{eq:DeltaweakIV}
\end{equation}

Generalizing from the example $\phi=1/t$, consider the sequence
\begin{align*}
r_n &= \sqrt{n}\\
\epsilon'_n &= \sqrt{n}\\
m_n &= \left (1+ \frac{1}{2\sqrt{n}}, \frac{1}{\sqrt{n}}\right )\\
A &= ]-\tfrac{1}{2}, \tfrac{1}{2}[ \;\times \; ]-2,-1[.
\end{align*}
For any $t\in m_n + \tfrac{1}{r_n} A$ we get
\[\Delta(t, m_n) \geq \sqrt{n} \cdot \frac{3}{2},\]
which implies that the condition of theorem \ref{theo:uniformdelta}, part 2 is fulfilled, and the delta method is not uniformly valid in this setting.

\subsection{Moment inequalities}

Suppose we observe an i.i.d. sample $(Y_{i1}, Y_{i2})$,
where we assume for simplicity that $\Var(Y)=I$, and $E[Y] \neq 0$.
We are interested in testing the hypothesis $E[Y] \geq 0$.\footnote{Throughout this section, inequalities for vectors are taken to hold for all components.}
A leading example of such a testing problem is inference under partial identification as discussed in \cite{imbens2004confidence}.
We follow the approach of \cite{hahnridder2013} in considering this as a joint hypothesis testing problem, and using a likelihood ratio test statistic based on the normal limit experiment. We demonstrate that the ``naive'' approximation to the distribution of this statistic, a $0.5\cdot \chi^2_1$ distribution is not uniformly valid, even though it \emph{is} pointwise valid.

The log generalized likelihood ratio test statistic (equivalently in this setting, the Wald test statistic) takes the form
\[n \cdot \min_{t' \in \mathbb{R}^{2+}} \|t'-T_n\|^2,\]
where $T_n$ is the sample average of $Y_i$.

Let
\begin{align}
T_n &= E_n[Y] \nonumber\\
\mu(\theta) &=  E[Y] \nonumber\\
\phi_1(t) &= \left (\argmin_{t' \in \mathbb{R}^{2+}} \|t'-t\|^2 \right ) - t \textrm{, and} \nonumber\\
\phi_2(t) &= \min_{t' \in \mathbb{R}^{2+}} \|t'-t\|^2 =  \|\phi_1(t)\|^2 . 
\end{align}
The parameter space $\Theta$ we consider is the space of all distributions of $Y$ such that  $E[Y_1]=0$ or  $E[Y_2]=0$, but  $E[Y]\neq 0$.
 
``Conventional'' testing of the null hypothesis $E[Y] \geq 0$ is based on critical values for the distribution of $n \cdot \phi_2(T_n)$, which are derived by applying (i) a version of the delta method to obtain the asymptotic distribution of $\phi_1(T_n)$, and (ii) the continuous mapping theorem to obtain the distribution of $\phi_2(T_n)$ from the distribution of $\phi_1(T_n)$.

We could show that the delta method does not yield uniformly valid approximations on $\Theta$ for  $\phi_1(T_n)$, using the condition of theorem \ref{theo:uniformdelta}.
Standard approximations in this setting, however, do not exactly fit into the framework of the delta method, since they \emph{do} account for the fact that one-sided testing creates a kink in the mapping $\phi_1$ at points where $t_1=0$ or $t_2=0$. 
As a consequence, it  is easier to explicitly calculate the remainder of the standard approximation and verify directly that this remainder does not vanish uniformly, so that uniform convergence does not hold for the implied approximations of $\phi_1(T_n)$ and $\phi_2(T_n)$.
We can rewrite
\begin{align*}
\phi_1(t) &= (\max(-t_1, 0), \max(-t_2, 0))\\
\phi_2(t) &= t_1^2 \cdot \mathbf{1}(t_1 \leq 0) + t_2^2 \cdot \mathbf{1}(t_2 \leq 0).
\end{align*}

Consider without loss of generality the case $m_2= E[Y_2] >0$ and $m_1= E[Y_1]=0$. For this case, the ``conventional'' asymptotic approximation\footnote{We can interpret this approximation as a first order approximation based on directional derivatives.} sets
\[\phi_1(t) \approx \tilde{\phi}_1(t):=   (\max(-t_1, 0), 0).\]
Based on this approximation, we obtain the pointwise asymptotic distribution of $n \cdot \phi_2(T_n)$ as $0.5 \chi^2_1 + 0.5 \delta_0$.
The remainders of these approximations are independent of the approximations themselves, and are equal to 
\begin{align*}
\phi_1(t) -  \tilde{\phi}_1(t) &= (0, \max(-t_2, 0))\\
\phi_2(t) - \tilde{\phi}_2(t) &= t_2^2 \cdot \mathbf{1}(t_2 \leq 0).
\end{align*}
These remainders, appropriately rescaled, do converge to $0$ in probability pointwise on $\Theta$, since $\sqrt{n} (T_2 - m_2) \rightarrow N(0,1)$.
This convergence is not uniform, however.
Consider a sequence of $\theta_n$ such that $m_{n1} =0$ and $m_{n2}=1/n$. For such a sequence we get
\begin{align*}
\phi_1(T_n) -  \tilde{\phi}_1(T_n) & \rightarrow^d (0, \max(Z,0)) \\
\phi_2(T_n) - \tilde{\phi}_2(T_n) &\rightarrow^d Z^2 \cdot \mathbf{1}(Z \leq 0).
\end{align*}
where
\[Z\sim N(1,1).\]

%
%

\subsection{Minimum distance estimation}

Suppose that we have (i) an estimator $T_n$ of various reduced-form moments $\mu(\theta)$, and that (ii) we also have a structural model which makes predictions about these reduced form moments $\mu(\theta)$.
If the true structural parameters are equal to $\beta$, then the reduced form moments are equal to $m(\beta)$.
Such structural models are often estimated using minimum distance estimation. Minimum distance finds the estimate $X_n$ of $\beta$ such that $m(X_n)$ gets as close as possible to the estimated moments $T_n$.

If the model is just-identified, we have $\dim(t) =\dim(x)$ and the mapping $m$ is invertible. In that case we can set
\[X_n = \phi(T_n) = m^{-1}(T_n),\]
and our general discussion applies immediately.

If the model is over-identified, there are more reduced form moments than structural parameters.
For simplicity and specificity, assume that there are two reduced form moments but only one structural parameter, $\dim(t) =2 > \dim(x) = 1$. 
Suppose that $T_n$ converges uniformly in distribution to a normal limit,
\begin{align*}
\sqrt{n} \cdot (T_n - \mu(\theta)) &\rightarrow^d N(0, \Sigma(\theta)) .
\end{align*}
Let $X_n$ be the (unweighted) minimum distance estimator of $\beta$, that is
\[
X_n = \phi(T_n) =\argmin_x e(x,T_n),
\]
where
\[e(x,T_n)= \|T_n -m(x)\|^2.\]

A delta-method approximation of the distribution of $X_n$ requires the slope of the mapping $\phi$. We get this slope by applying the implicit mapping theorem to the first-order condition $\partial_x e =0$.  This yields 
\begin{align*}
\partial_x e &= - 2 \cdot\partial_x m \cdot (t - m(x))\\
\partial_{xx}  e &=	-2 \cdot(\partial_{xx}  m \cdot (t-m) - \|\partial_x m\|^2)	\\
\partial_{x t} e &=	- 2 \cdot\partial_x m	\\
\partial_t \phi (t) &= - (\partial_{xx}  e)^{-1} \cdot \partial_{x t} e\\
&= -  (\partial_{xx}  m \cdot (t-m) - \|\partial_x m\|^2)^{-1}  \cdot \partial_x m. 
\end{align*}

If the model is correctly specified, then there exists a parameter value $x$ such that $m(x) = \mu(\theta)$.
Evaluating the derivative $\partial_t \phi (t)$ at $t= m(x)$ yields
\begin{align*} 
\partial_t \phi (t) 
&=    \frac{1}{\|\partial_x m\|^2}   \cdot \partial_x m. 
\end{align*}
Let $m=m(x)$, so that $\phi(m) = x$.
The normalized remainder of a first order approximation to $\phi$ at such a point is given by
\begin{align}
\Delta(t, m) &= \frac{1}{\|\partial_m\phi(m)\| \cdot\|t-m\|}\left |\phi(t)- \phi(m) -  \partial_m\phi(m) \cdot(t - m) \right | \nonumber \\
&=   \frac{\|\partial_x m\|^{2}}{\|\partial_x m\| \cdot\|t-m\|}\left |\phi(t)- x -  \|\partial_x m\|^{-2}  \cdot \partial_x m \cdot(t - m) \right |.
\end{align}

The magnitude of this remainder depends on the curvature of the manifold traced by $m(.)$, as well as on the parametrization which maps $x$ to this manifold.
The remainder will be non-negligible to the extent that either the manifold or the parametrization deviate from linearity.
If the manifold has kinks, that is points of non-differentiability, then our discussion of moment inequalities immediately applies.
If the manifold is smooth but has a high curvature, then the delta-method will provide poor approximations in finite samples, as well.
As a practical approach, we suggest to numerically evaluate $\Delta$ for a range of plausible values of $m$ and $t$.

%
%
%
%
%

\section{Conclusion}

 Questions regarding the uniform validity of statistical procedures figure prominently in the econometrics literature in recent years: Many conventional procedures perform poorly for some parameter configurations, for any given sample size, despite being asymptotically valid for all parameter values.
 We argue that a central reason for such lack of uniform validity of asymptotic approximations rests in failures of the delta-method to be uniformly valid. 
 
 In this paper, we provide a condition which is both necessary and sufficient for uniform negligibility of the remainder of delta-method type approximations. This condition involves a uniform bound on the behavior of the remainder of a Taylor approximation.
 We demonstrate in a number of examples that this condition is fairly straightforward to check, either analytically or numerically. The stylized examples we consider, and for which our necessary condition fails to hold, include $1/t$, $t^2$,  $|t|$, and $\sqrt{t}$. In each of these cases problems arise in a neighborhood of $t=0$. Problems can also arise for large $t$.
We finally discuss three more realistic examples, weak instruments, moment inequalities, and minimum distance estimation.

\clearpage

\appendix

\section{Proofs}
\label{sec:proofs}

\textbf{Proof of lemma \ref{lem:characterization}:}
\begin{enumerate}
\item 
To see that convergence along any sequence $\theta_n$ follows from this condition, note that
\[ \sup_\theta d^{\theta}_{BL} (X_n, Y_n) 
\geq   d^{\theta_n}_{BL} (X_n, Y_n) .\]

To see that  convergence along any sequence implies this condition, note that 
\[ \sup_\theta d^{\theta}_{BL} (X_n, Y_n)  \not\rightarrow 0\]
implies that there exist $\epsilon > 0$, and sequences $\theta_m$, $n_m \rightarrow \infty$, such that
\[ d^{\theta_m}_{BL} (X_{n_m}, Y_{n_m})>\epsilon\]
for all $m$.

\item Similarly 
\[\sup_\theta P^\theta(\|X_n -Y_n\| >\epsilon) \geq P^{\theta_n}(\|X_n - Y_n\| >\epsilon) \]
shows sufficiency, and
\[\sup_\theta P^\theta(\|X_n -Y_n\| >\epsilon)  \not\rightarrow 0 \]
implies that there exist $\epsilon'>0$ and sequences $\theta_m$, $n_m$, such that
\[P^{\theta_m}(\|X_{n_m} - Y_{n_m}\| >\epsilon)>\epsilon' \]
for all $m$.
\end{enumerate}
$\Box$\\[10pt]


 \textbf{Proof of lemma \ref{lem:uniforminference}: }
 Fix an arbitrary sequence $\theta_n$. Uniform convergence in distribution of $Z_n$ to $Z$ implies convergence in distribution of $Z_n$ to $Z$ along this sequence.
 By Portmanteau's lemma (\citealt{van2000asymptotic}, p6), uniform convergence in distribution of $Z_n$ to $Z$ is equivalent to convergence of $F^{\theta_n}_{Z_n}(\widehat{z})$ to $F_Z(\widehat{z})$ at all continuity points $\widehat{z}$ of $F_Z(.)$. Since we assume the latter to be continuous, convergence holds at all points $\widehat{z}$, and thus in particular at the critical value $z$. The claim follows.
$\Box$\\

\textbf{Proof of theorem \ref{theo:uniformCMT}:}\\
Let $1\leq \kappa <\infty$ be such that  $|\psi(x) - \psi(y)| \leq \kappa \cdot \|x-y\|$ for all $x,y$.
\begin{enumerate}
\item Note that $h\in BL_1$ implies $h':=\tfrac{1}{\kappa} \cdot h\circ \psi \in BL_1$:
\begin{align*}
| h'(x) - h'(y) | &= \kappa^{-1} | h(\psi(x))-h(\psi(y))| \\
& \text{(by definition of $h'$)}\\
& \leq \kappa^{-1} ||\psi(x)-\psi(y)|| \\
& \text{(since $h \in BL_1$)} \\ 
& \leq  \kappa^{-1} \kappa || x-y || \\
& \text{(since $\psi$ is Lipschitz-continuous with parameter $\kappa$),}
\end{align*}
and $|h'(x)| \leq 1$ for $\kappa \geq 1$.

By definition of the Lipschitz metric
\begin{align*}
d^{\theta}_{BL} (\psi(X_n), \psi(Y_n)) &=\sup_{h \in \mathbf{BL_1}} \left | E^{\theta}[h(\psi(X_n))] - E^{\theta}[h(\psi(Y_n))] \right |\\
&\leq \kappa  \cdot \sup_{h' \in \mathbf{BL_1}} \left | E^{\theta}[h'(X_n)] - E^{\theta}[h'(Y_n)] \right |\\
&= \kappa \cdot d^{\theta}_{BL} (X_n, Y_n). 
\end{align*}
$d^{\theta_n}_{BL} (X_n, Y_n)  \rightarrow 0 $ for all sequences $\{\theta_n \in \Theta\}$ therefore implies
$d^{\theta_n}_{BL} (\psi(X_n), \psi(Y_n))  \rightarrow 0 $ for all such sequences.

\item For a given $\epsilon > 0$, let $\delta >0$ be such that
$\|x-y\| \leq \delta$ implies $\|\psi(x) - \psi(y) \| \leq \epsilon$ for all $x,y$.
Such a $\delta$ exists by uniform continuity of $\psi$.
For this choice of $\delta$, we get
\[P^\theta(\|\psi(X_n) - \psi(Y_n)\| >\epsilon)  \leq P^\theta(\|X_n - Y_n\| > \delta ).  \]
By uniform convergence in probability,
 $P^{\theta_n}(\|X_n - Y_n\| > \delta) \rightarrow 0$
for all $\delta > 0$ and all sequences $\{\theta_n \in \Theta\}$,  which implies
 $P^{\theta_n}(\|\psi(X_n) - \psi(Y_n)\| >\epsilon) \rightarrow 0 $ for all such sequences.
\end{enumerate}
$\Box$\\[10pt]

\textbf{Proof of theorem \ref{theo:uniformdelta}:}\\   
Define 
\begin{align*}
\widetilde{X}_n &= E(\mu) D(\mu) \cdot S_n \textrm{ and} \\
R_n &=X_n - \widetilde{X}_n.
\end{align*} 
The proof is structured as follows.
We show first that under our assumptions $\widetilde{X}_n$ converges uniformly in distribution to $X = E(\mu) D(\mu) \cdot S$. This is a consequence of uniform convergence in distribution of $S_n$ to $S$ and the boundedness of $E(\mu) D(\mu)$.

We then show that $R_n$ converges to $0$ in probability uniformly under the sufficient condition \eqref{eq:vanishingremainder}. This, in combination with the first result, implies uniform convergence in distribution of $X_n$ to $X$, by Slutsky's theorem applied along arbitrary sequences $\theta_n$.

We finally show that $R_n$ diverges along some sequence $\theta_n$ under condition \eqref{eq:divergingremainder}. This implies that $X_n = \widetilde{X}_n + R_n$ cannot converge in distribution along this sequence.\\

\textbf{Uniform convergence in distribution of $\widetilde{X}_n$ to $X$:}\\
Note that 
\[d^\theta_{BL}(\widetilde{X}_n, X) \leq d_x\cdot d^\theta_{BL}(S_n, S).\]
This holds since multiplication by $E(\mu) D(\mu)$ is a Lipschitz continuous function with Lipschitz constant $d_x$. To see this note that each of the $d_x$ rows of $E(\mu) D(\mu)$ has norm $1$ by construction.

Since $d^\theta_{BL}(S_n, S) \rightarrow 0$ by assumption, the same holds for $d^\theta_{BL}(\widetilde{X}_n, X) $.\\

\textbf{Uniform convergence in probability of $X_n$ to $\widetilde{X}_n$ under condition \eqref{eq:vanishingremainder}:}\\
We can write
\begin{align}
\|R_n \| &= \| \Delta(T_n, \mu) \| \cdot r_n \|T_n - \mu\| \nonumber\\
&=  \|  \Delta(\mu + S_n/r_n, \mu) \| \cdot  \|S_n\| \label{eq:remainder}\\
&\leq \delta(\|S_n\|/r_n) \cdot \|S_n\|. \nonumber
\end{align}
Fix $M$ such that $P_\theta(\|S\| >M) < \epsilon / 2$ for all $\theta$; this is possible by tightness of $S$ as imposed in assumption \ref{as:uniformS}.
By uniform convergence in distribution of $S_n$ to to the continuously distributed $S$, this implies $P_{\theta_n}(\|S_n\| >M) < \epsilon$ for any sequence $\theta_n \subset \Theta$ and $n$ large enough.
We get
\begin{align*}
P_{\theta_n}(\|R_n \|> \epsilon) &\leq P_{\theta_n}(\|S_n\| > M) + P_{\theta_n}(\|S_n\| \leq  M, \;\delta(M /r_n) >\epsilon/ M)\\
&=  P_{\theta_n}(\|S_n\| > M) <\epsilon.
\end{align*}
for any sequence $\theta_n \subset \Theta$ and $n$ large enough, using condition \eqref{eq:vanishingremainder}. But this implies $P_{\theta_n}(\|R_n \|> \epsilon) \rightarrow 0$. \\

\textbf{Existence of a diverging sequence under condition \eqref{eq:divergingremainder}:}\\
Let $\theta_n$ be such that $\mu(\theta_n) = m_n$, where $(m_n, \epsilon'_n)$ is a sequence such  that condition \eqref{eq:divergingremainder} holds.
By equation \eqref{eq:remainder},
\begin{align*}
P_{\theta_n}(\|R_n \|> \epsilon'_n / \underline{s}) &\geq P_{\theta_n}(S_n \in A, \; \Delta(m_n + S_n/r_n, m_n) > \epsilon_n')\\
& = P_{\theta_n}(S_n \in A) > \underline{p}/2
\end{align*}
for $n$ large enough, under the conditions imposed, using again uniform convergence in distribution of $S_n$ to $S$.\\

Note that $R_n = X_n - \widetilde{X}_n$ and thus $\|R_n\| \leq \|X_n\|  + \|\widetilde{X}_n\|$,
which implies
\[P( \|R_n\|  > \epsilon'_n / \underline{s}) \leq P( \|X_n\|  > \epsilon'_n / (2 \underline{s})) +  P( \|\widetilde{X}_n\|  > \epsilon'_n / (2 \underline{s})).\]
Suppose that $X_n \rightarrow^d X$ (and recall $\widetilde{X}_n \rightarrow^d X$) for the given sequence $\theta_n$.
Since $X$ is tight and $\epsilon'_n \rightarrow \infty$, this implies $P( \|X_n\|  > \epsilon'_n / (2 \underline{s})) \rightarrow 0$, similarly for $\widetilde{X}_n$, and thus
\[P( \|R_n\|  > \epsilon'_n / \underline{s}) \rightarrow 0.\] 
Contradiction. This implies that we cannot have $X_n \rightarrow^d X$ for the given sequence $\theta_n$.
$\Box$\\

\textbf{Proof of theorem \ref{theo:sufficient}: }\\
If $\mu(\Theta) \subset \mathbb{R}^l$ is compact, then so is $\mu(\Theta)^\epsilon$.
Since $\partial_t \phi $ is assumed to be continuous on $\mu(\Theta)^\epsilon$, it follows immediately that $\partial_t \phi $ is bounded on this set, and we also get $\|E(\mu)\| \leq \overline{E}$ for all $\mu \in \mu(\Theta)^\epsilon$.
Theorem 4.19 in \cite{rudin1964principles} furthermore implies that $\partial_t \phi $ is uniformly continuous on $\mu(\Theta)^\epsilon$.\\

Consider now first the case $\dim(x) = 1$.
Suppose $\|t_1 - m\| \leq \epsilon$ and $t_1, m \in \mu(\Theta)$.
By continuous differentiability and the mean value theorem, we can write 
\[\phi(t_1)- \phi(m) = \partial_t \phi (t_2) \cdot(t_1 - m)\]
 where 
 \[t_2 = \alpha t_1 + (1-\alpha) m \in \mu(\Theta)^\epsilon\]
 and $\alpha \in [0,1]$.
  We  get
\begin{align*}
\Delta(t_1, m) &= \frac{1}{\|m - t_1\|}\left \|E(m)\cdot(\phi(t_1)- \phi(m) -  \partial_m \phi (m) \cdot(t_1 - m)) \right \|\\
 &\leq  \frac{\overline{E}}{\|m - t_1\|}\left \|\left (\partial_t \phi(t_2) - \partial_m \phi(m)\right ) \cdot(t_1 - m) \right \|\\
 &\leq \overline{E} \cdot \left \|\partial_t \phi(t_2) -\partial_m \phi(m)\right \|.
\end{align*}
Uniform continuity of $\partial_{t} \phi$ implies that for  any $\delta>0$ there is an $\epsilon'$ such that
 $||t_2-m|| < \epsilon'$ implies that $\overline{E} \cdot|\partial_t \phi(t_2) -\partial_m \phi(m)|< \delta.$  
  Since $\|m - t_2\| \leq \|m - t_1\|$,
this implies that there exists a function $\delta(.)$ such that
\[\Delta(t_1, m) \leq \delta(\|m - t_1\|),\]
 and  $\delta$ goes to $0$ as its argument goes to $0$,
so that condition \eqref{eq:vanishingremainder} is satisfied.\\

Let us now consider the case $\dim(x) = d_x > 1$.
By the same arguments as for the case $\dim(x) = 1$, we get for the $i$th component of $\phi$ that
\begin{align*}
\Delta_i(t_1, m) &:= 
\frac{E_i}{\|m - t_1\|}\left \|\phi_i(t_1)- \phi_i(m) -  \partial_m \phi_i (m) \cdot(t_1 - m) \right \|\\
 &\leq  \overline{E} \cdot\left \|\partial_t \phi_i(t_{2,i}) -\partial_m \phi(m)\right \|.
\end{align*}
where 
 \[t_{2,i} = \alpha_i t_1 + (1-\alpha_i) m \in \mu(\Theta)^\epsilon.\]
As before, uniform continuity of $\partial_t \phi_i$ implies existence of a function $\delta_i$ such that
\[\Delta_i(t_1, m) \leq \delta_i (\|m - t_1\|),\]
 and  $\delta_i$ goes to $0$ as its argument goes to $0$.
By construction
\[\Delta(t_1, m) = \sqrt{\sum_i \Delta_i(t_1, m)^2} \leq d_x \cdot \max_i \Delta_i(t_1, m). \]
Setting
\[\delta(.) = d_x\cdot \max_i \delta_i(.)\]
then yields a function $\delta(.)$ which satisfies the required condition \eqref{eq:vanishingremainder}.
$\Box$\\

\textbf{Proof of lemma \ref{lem:uniformCLT}: }\\
By definition \ref{def:uniform}, we need to show convergence in distribution (ie., convergence with respect to the bounded Lipschitz metric) along arbitrary sequences $\theta_n$.

Consider such a sequence $\theta_n$, and define $\widetilde{Y}_{in} := \Sigma^{-1/2}(\theta_n) \cdot (Y_i  - \mu(\theta_n))$, so that $\Var(\widetilde{Y}_{in}) = I$.
Then the triangular array $\{\widetilde{Y}_{1n}, \ldots, \widetilde{Y}_{nn}\}$ with distribution corresponding to $\theta_n$ satisfies the conditions of Lyapunov's central limit theorem, and thus those of the Lindeberg-Feller CLT as stated in \cite[][proposition 2.27, p20]{van2000asymptotic}.
 We therefore have 
\[\widetilde{S}_n:= \frac{1}{\sqrt{n}} \sum_{i=1}^n \widetilde{Y}_{in} \rightarrow \widetilde{Z}\sim N(0,I)\]
in distribution, and thus with respect to the bounded Lipschitz metric, that is
\[d^{\theta_n}_{BL}( \widetilde{S}_n, \widetilde{Z} ) \rightarrow 0.\]

Now consider $S_n =\Sigma^{1/2}(\theta_n) \cdot  \widetilde{S}_n$
and $Z = \Sigma^{1/2}(\theta_n) \cdot  \widetilde{Z}$. 
We have
\[d^{\theta}_{BL}( S_n, Z) \leq \| \Sigma^{1/2}(\theta_n) \|   \cdot  d^{\theta}_{BL}( \widetilde{S}_n,  \widetilde{Z})\]
-- this follows from the definition of the bounded Lipschitz metric, again by the same argument as in the proof of theorem \ref{theo:uniformCMT}.
Since $d^{\theta}_{BL}( \widetilde{S}_n, \widetilde{Z} ) \rightarrow 0$, and  $\| \Sigma^{1/2} \|$ is bounded by assumption, we get $d^{\theta}_{BL}( S_n, Z) \rightarrow 0$, and the claim follows. $\Box$\\

\textbf{Proof of proposition \ref{prop:normalpivot}: }\\
This follows immediately from theorem \ref{theo:uniformdelta} and Slutsky's theorem, applied along any sequence of $\theta_n$.
$\Box$\\

\textbf{Derivation of equation \eqref{eq:DeltaweakIV}: } 
\begin{eqnarray*}
&\quad & \frac{|| t-m ||}{E(m)} \Delta(t,m)= \\
&=& \Big| \frac{t_1}{t_2} - \frac{m_1}{m_2} - \frac{1}{m_2} (t_1-m_1) + \frac{m_1}{m^2_2} (t_2-m_2) \Big| \\
&=& \Big| \frac{t_1m_2-t_2m_1}{ t_2 m_2} - \frac{1}{m_2} (t_1-m_1) + \frac{m_1}{m^2_2} (t_2-m_2) \Big| \\
&=& \Big| \frac{t_1m_2-m_1m_2+m_1m_2-t_2m_1}{ t_2 m_2} - \frac{1}{m_2} (t_1-m_1) + \frac{m_1}{m^2_2} (t_2-m_2) \Big| \\
&= &\Big| \frac{1}{ t_2}(t_1-m_1) -\frac{m_1}{t_2 m_2} (t_2-m_2)- \frac{1}{m_2} (t_1-m_1) + \frac{m_1}{m^2_2} (t_2-m_2) \Big| \\
&=& \Big| (t_1-m_1) \Big[ \frac{1}{t_2}-\frac{1}{m_2} \Big]  - \frac{m_1}{m_2} (t_2-m_2) \Big[ \frac{1}{t_2}-\frac{1}{m_2} \Big] \Big| \\
&=& \Big| \frac{t_2-m_2}{t_2 m_2} \Big| \: \: \Big| (t_1-m_1) - \frac{m_1}{m_2}(t_2-m_2) \Big|.
\end{eqnarray*}
$\Box$\\

\clearpage 
\bibliographystyle{apalike}
\bibliography{Uniformity}

\end{document}